\documentclass[12pt]{article}
\usepackage{amsmath,amssymb}
\newcommand{\qed}{{\unskip\nobreak\hfil\penalty50\hskip2em\vadjust{}
    \nobreak\hfil$\square$\parfillskip=0pt\finalhyphendemerits=0\par}}

\newtheorem{thm}{Theorem}
\newcommand{\nn}{_{\langle n\rangle}}

\newcommand{\<}{\langle}
\renewcommand{\>}{\rangle}

\newcommand{\bs}{\bigskip}

\newcommand{\C}{\mathbb C}

\renewcommand{\d}{\delta}

\newcommand{\J}{\mathcal J}

\newcommand{\mb}{\mbox}

\newcommand{\ot}{\otimes}
\newcommand{\p}{\partial}

\newcommand{\stk}{\stackrel}
\renewcommand{\th}{\theta}

\newcommand{{\z}}{\Bbb Z}

\begin{document}
\setlength{\baselineskip}{18pt}

\begin{center}
{\Large\bf A Note on Cyclic Gradients}\bs\bs

{\sc Dan Voiculescu}

\bs\bs\bs {\it To the memory of Gian-Carlo Rota}
\end{center}

\bs\bs
The cyclic derivative was introduced by G.-C.~Rota, B.~Sagan and 
P.~R.~Stein in [3]
as an extension of the derivative to noncommutative polynomials. Here 
we show that
there are simple necessary and sufficient conditions for  an 
$n$-tuple of polynomials
in $n$ noncommuting indeterminates to be a cyclic gradient (see Theorem 1) and
similarly for a polynomial to have vanishing cyclic gradient (see Theorem 2).
  Our interest in cyclic
gradients stems from free probability theory and random matrices (see 
the Remark at
the end) [1],[2],[4],[5],[6]. This note should also reduce the
paucity of results on cyclic derivatives in several variables pointed 
out in [3,
page 73].

\bs
Let \ $K\nn = K\<X_1,\dots ,X_n\>$ \ be the ring of polynomials in
noncommuting indeterminates $X_1,\dots ,X_n$ with coefficients in the 
field $K$ of
characteristic zero. The partial generalized difference quotients are the
derivations
$$
\p_j: K\nn\to K\nn\ot K\nn
$$
such that $\p_j X_k=0$ \ if \ $j\neq k$ and $\p_jX_j=1\ot 1$.
The $\ot$ here is over $K$ and $K\nn\ot K\nn$ is given the bimodule 
structure such
that $a(b\ot c)=ab\ot c$, \ $(b\ot c)d=b\ot cd$.

The partial cyclic derivatives are then
$$
\d_j={\tilde\mu}\circ\p_j:K\nn\to K\nn
$$
where ${\tilde\mu}(a\ot b)=ba$.

We shall denote by $N:K\nn\to K\nn$ the ``number operator", i.e.~the
linear map so that $N1=0$, \
  $NX_{i_1}\dots X_{i_k}=kX_{i_1}\dots X_{i_k}$. Also,
$CK\nn$ will denote the cyclic subspace, i.e.~the vector subspace 
spanned by all
cyclic symmetrizations of monomials
$$
CX_{i_1}\dots X_{i_p} = \sum_{1\leq j\leq p}
X_{i_{j+1}}\dots X_{i_p}X_{i_1}\dots X_{i_j} \ , \qquad
p\!\geq\! 1 \ \ {\mb{and}} \ \ C1\!=\!0
$$
(the constants are not in the cyclic subspace).

\newpage
\begin{thm} Let $P_1,\dots ,P_n\in K\nn$. The following conditions are
equivalent:

{\rm{(i)}} \ there is $P\in K\nn$ such that $\d_jP=P_j$ $(1\leq j\leq n)$.

{\rm{(ii)}} \ $\displaystyle{\sum_{1\leq j\leq n}} [X_j,P_j]=0$.

{\rm{(iii)}} \ $\displaystyle{\sum_{1\leq j\leq n}} X_jP_j\in CK\nn$.

{\rm{(iv)}} \ $\d_k\displaystyle{\left(\sum_{1\leq j\leq n} X_jP_j\right)}
=(N+I)P_k$. (Here $I$ denotes the identity map of $K\nn$ to itself.)
\end{thm}

{\bf Proof.} It is easily seen that it suffices to prove the theorem 
for homogeneous
$P_1,\dots ,P_n$ of the same degree, i.e.~we may assume $NP_j=sP_j$ 
$(1\leq j\leq
n)$ for some $s\geq 0$. Also the case of constants being obvious we 
will concentrate
on $s\geq 1$.

(i)\,$\Rightarrow$\,(ii) To check that
$\displaystyle{\sum_{1\leq j\leq n}} [X_j,\d_jP]=0$ \ it suffices to 
do so when $P$
is a monomial $X_{i_0}\dots X_{i_s}$. Then
$$
\d_j P=\sum_{i_p=j} X_{i_{p+1}}\dots X_{i_s}X_{i_0}\dots X_{i_{p-1}}
$$
so that
$$
[X_j,\d_jP] = \sum_{i_p=j}
\big(X_{i_p}\dots X_{i_s}X_{i_0}\dots  X_{i_{p-1}}-
  X_{i_{p+1}} \dots X_{i_s}X_{i_0}\dots  X_{i_p}\big)
$$
and hence
$$
\sum_{1\leq j\leq n} [X_j,\d_jP]
= \sum_{1\leq p\leq s}
\big(X_{i_p}\dots X_{i_s}X_{i_0}\dots  X_{i_{p-1}}-
  X_{i_{p+1}} \dots X_{i_s}X_{i_0}\dots  X_{i_p}\big) =0 \ .
$$

(ii)\,$\Rightarrow$\,(iii) \ Let
$$
P_j=\sum_{i_1,\dots ,i_s} c^j_{i_1\dots i_s} X_{i_1}\dots  X_{i_s}
$$
The coefficient of $X_{i_0}\dots  X_{i_s}$ in
$\displaystyle{\sum_{1\leq j\leq n}} [X_j,P_j]$  \ is
$$
c^{i_0}_{i_1\dots i_s} - c^{i_s}_{i_0\dots i_{s-1}} \ .
$$
Hence (ii) gives $c^{i_0}_{i_1\dots i_s} = c^{i_s}_{i_0\dots i_{s-1}}$.
On the other hand, if $c_{i_0\dots i_s}$ denotes the coefficient of
$X_{i_0}\dots X_{i_s}$ in $\displaystyle{\sum_{i\leq j\leq n}}X_jP_j$, clearly
$c_{i_0\dots i_s} = c^{i_0}_{i_1\dots i_s}$, so that (ii) implies
$c_{i_0\dots i_s} = c_{i_si_0\dots i_{s-1}}$, i.e.~cyclicity.

(iii)\,$\Rightarrow$\,(iv) As before, let $c^j_{i_1\dots i_s}$ and
$c_{i_0\dots i_s}$ denote the coefficients of $P_j$ and $\displaystyle{\sum_j}
X_jP_j$  respectively. Then $c^{i_0}_{i_1\dots i_s} = c_{i_0\dots i_s}$
and the cyclicity condition gives
$c_{i_0\dots i_s} = c_{i_si_0\dots i_{s-1}}$. We have
\begin{eqnarray*}
&& \d_k\left(\sum_k X_jP_j\right) \\
&& \quad = \sum_{i_0\dots i_s} \sum_{\{r:i_r=k\}}
c_{i_0\dots i_s} X_{i_{r+1}}\dots X_{i_s}X_{i_0}\dots X_{i_{r-1}} \\
&& \quad =  \sum_{i_0\dots i_s} \sum_{\{r:i_r=k\}}
c^k_{i_{r+1}\dots i_s i_0\dots i_{r-1}} X_{i_{r+1}}\dots X_{i_s}X_{i_0}\dots
X_{i_{r-1}} \\
&& \quad =
  \sum_{0\leq r\leq s} \ \sum_{i_0\dots i_{r-1}i_{r+1}\dots i_s}
c^k_{i_{r+1}\dots i_s i_0\dots i_{r-1}}
  X_{i_{r+1}}\dots X_{i_s}X_{i_0}\dots X_{i_{r-1}}
= (s+1)P_k \ .
\end{eqnarray*}

(iv)\,$\Rightarrow$\,(i)
Since the $P_j$ are homogeneous of the same degree and the field 
characteristic is
zero, this is obvious. \qed

\bs
There is also a simple description of the noncommutative polynomials with
vanishing cyclic gradient.

\begin{thm} We have
$$
{\rm{Ker}} \ \d \ =\sum_{1\leq k\leq n} [X_k,K\nn]+{\C}1  \
= \ {\C}1 + [K\nn,K\nn] \ =  \ {\rm{Ker}} \ C
$$
\end{thm}

{\bf Proof.} (i) Ker $\d\subset {\mb{Ker}} \ C$. We have
$$
Cp \ =\sum_{1\leq j\leq m} X_j\d_j p \ = \ 0
$$

(ii) Clearly,
$$
\sum_{1\leq k\leq n} [X_k,K\nn]+{\C}1 \subset [K\nn,K\nn] +{\C}1
$$
Also, since $1\in {\mb{Ker}} \ C$ and
$[X_{i_1}\dots X_{i_r},X_{i_{r+1}}\dots X_{i_{r+s}}]$ is the difference
of two cyclic permutations of $X_{i_1}\dots X_{i_{r+s}}$, we have
$[K\nn,K\nn]+{\C}1\subset{\mb{Ker}} \ C$.

To see that \  Ker $C\subset\sum_{1\leq k\leq n} [X_k,K\nn]+{\C}1$,
remark that \ Ker $C$ is spanned by homogeneous elements and that
$Cp=0$, where $p$ is homogeneous of degree $m$ iff $p$ is a linear
combination of differences
$X_{i_1}\dots X_{i_m}-X_{i_2}\dots X_{i_m}X_{i_0}$.

(iii) To see that
$\displaystyle{\sum_{1\leq k\leq n}} [X_k,K\nn]+ \ {\C}1\subset 
{\mb{Ker}} \ \d$,
it suffices to show that\\$[X_k,X_{i_1}\dots X_{i_s}]\in {\mb{Ker}} \ \d$.
This is clearly so, since
$$
[X_k,X_{i_1}\dots X_{i_s}] = X_k X_{i_s}\dots X_{i_s}-X_{i_1}\dots X_{i_s}X_k
$$
and the cyclically equivalent elements
$X_k X_{i_1}\dots X_{i_s}$,  $X_{i_1}\dots X_{i_s}X_{i_k}$ have the
same cyclic gradient. \qed

\bs
Putting together the two theorems, we have an exact sequence
$$
0\to [K\nn,K\nn]\to K\nn \stk{\d}{\longrightarrow}
(K\nn)^n \stk{\th}{\longrightarrow} K\nn
$$
where \ $\th((P_j)_{1\leq j\leq n})=\displaystyle{\sum_j}\, [X_j,P_j]$.

\bs

{\bf Remark.} The motivation for this note is from free entropy and 
large deviations
for random matrices. Let $(M,\tau)$ be a von Neumann algebra with 
normal faithful
trace-state $\tau$ and $X_k=X^*_k\in M$ \ $(1\leq k\leq n)$ which are 
algebraically
free.

\bs
Let ${\J}_k={\J}(X_k:{\C}\< X_1,\dots ,X_{k-1},X_{k+1},\dots ,X_n\>)$
be the noncommutative Hilbert transforms defined in [4], in 
connection with free
entropy. On the other hand, the upper bound for large deviations for 
$n$-tuples of
random matrices found in [2] fits well with free entropy except for a 
term involving
cyclic gradients and about which it is not known whether it is not 
actually zero.
The precise question is, whether the $n$-tuple $({\J}_k)_{1\leq k\leq n}$
(when it exists) is a limit in 2-norm of cyclic gradients of polynomials in
the noncommuting variables $X_1,\dots ,X_n$\,? The theorem we proved 
here provides a
partial affirmative answer:
\begin{quote}
If \ $({\J}_k)_{1\leq k\leq n}$ \ are noncommutative polynomials in
$X_1,\dots ,X_n$, then there is a noncommutative polynomial $P$ in
$X_1,\dots, X_n$ such that ${\J}_k=\d_k P$ \ $(1\leq k\leq n)$.
\end{quote}
Indeed, by Corollary 5.12 in [5] we have $\displaystyle{\sum_k}\, 
[{\J}_k,X_k]=0$.
Hence the commutator condition (ii) in the Theorem is satisfied.

\bs {\bf Acknowledgments.} This research was conducted by the author for
the Clay Mathematics Institute. Partial support was also provided by
National Science Foundation Grant DMS95--00308.

\bs\begin{center}{\bf References}\end{center}
\begin{description}
\item{[1]} \ P.Biane, R.Speicher. Free diffusions, free entropy
and free Fisher information, preprint.
\item{[2]} \ T.Cabanal-Duvillard, A.Guionnet. Large deviations upper bounds
and noncommutative entropies for some matrices ensembles, preprint.
\item{[3]} \  G.-C.Rota, B.Sagan, P.R.Stein. A cyclic derivative in 
noncommutative
algebra.\\{\it Journal of Algebra} {\bf 64}, 54--75 (1980).
\item{[4]} \ D.Voiculescu. The analogues of entropy and of Fisher's
information measure in free probabilitiy theory, V: noncommutative Hilbert
transforms. {\it Invent. Math.} {\bf 32}, no.~1, 189--227 (1998).
\item{[5]} \  D.Voiculescu. The analogues of entropy and of Fisher's
information measure in free probabilitiy theory, VI: liberation
and mutual free information. {\it Advances in Mathematics} {\bf 146}, 101--166
(1999).
\item{[6]} \ D.Voiculescu. Lectures on free probability theory. Notes 
for a course
at the Saint-Flour Summer School on Probability Theory, preprint (1998).
\end{description}

\bs\bs
\mbox{}\hfill
\begin{tabular}{r}
Department of Mathematics\\University of California\\
Berkeley, California 94720--3840\\
{\tt dvv@math.berkeley.edu}
\end{tabular}

\end{document}